\documentclass[11pt]{article}
\usepackage{amsmath,amsthm,amssymb,amsfonts}
\usepackage[numbers,sort&compress]{natbib}
\usepackage{color,colordvi}
\usepackage{fullpage}

\newtheorem{theorem}{Theorem}[section]
\newtheorem{lemma}[theorem]{Lemma}%[section]
\newtheorem{corollary}[theorem]{Corollary}%[section]
\newtheorem{remark}{Remark}%[section]
\addtocounter{section}{0}
\begin{document}
%Spectral radius and factors of bipartite graphs
\title{\LARGE{\bf The existence of biregular spanning subgraphs in bipartite graphs via spectral radius}}
\author{Dandan Fan$^{a}$, Xiaofeng Gu$^b$, Huiqiu Lin$^{c}$
%\thanks{Corresponding author. \newline {\it E-mail}: {\tt ddfan0526@163.com} (D. Fan), {\tt xgu@westga.edu} (X. Gu), {\tt huiqiulin@126.com} (H. Lin)}
\\[2mm]
\small\it $^a$ College of Mathematics and Physics, Xinjiang Agricultural University,\\ \small\it Urumqi, Xinjiang 830052, China\\[1mm]
\small\it $^b$ College of Mathematics, Computing, and Sciences, University of West Georgia, \\ \small\it Carrollton, GA 30118, USA\\ [1mm]
\small\it $^c$ School of Mathematics, East China University of Science and Technology, \\ \small\it   Shanghai 200237, China\\\footnotesize {\it Email}: {\tt ddfan0526@163.com} (D. Fan), {\tt xgu@westga.edu} (X. Gu), {\tt huiqiulin@126.com} (H. Lin)
}
\date{}
\maketitle

\begin{abstract}
Biregular bipartite graphs have been proven to have similar edge distributions to random bipartite graphs and thus have nice pseudorandomness and expansion properties. Thus it is quite desirable to find a biregular bipartite spanning subgraph in a given bipartite graph. In fact, a theorem of Ore implies a structural characterization of such subgraphs in bipartite graphs. In this paper, we demonstrate the existence of biregular bipartite spanning subgraphs in bipartite graphs by employing spectral radius. We also study the existence of spanning trees with restricted degrees and edge-disjoint spanning trees in bipartite graphs via spectral radius. 

\end{abstract}


\section{Introduction}

An $(X, Y)$-bipartite graph is a bipartite graph with bipartition $(X, Y)$. An $(X, Y)$-bipartite graph is called {\it $(a,b)$-biregular} if each vertex in $X$ has degree $a$ and each vertex in $Y$ has degree $b$. An expander mixing lemma for biregular bipartite graphs was provided by Haemers in \cite[Theorem~3.1.1]{Haem79} and also in \cite[Theorem 5.1]{Haem95}. By this bipartite expander mixing lemma, biregular bipartite graphs have similar edge distributions to those of random bipartite graphs and thus are called {\it pseudorandom bipartite graphs}. Spectral expansion properties of these kinds of graphs have been studied in \cite{FGL24}.  Biregular bipartite graphs are aslo closely related to unbalanced bipartite Ramanujan graphs and unbalanced bipartite expander graphs (see \cite{EFMP23} for mroe details). Thus it is quite interesting and desirable to study the existence of biregular bipartite spanning subgraphs in a given bipartite graph.

In general, an \textit{$[a,b]$-factor} of a graph $G$ is a spanning subgraph $H$ such that $a\leq d_{H}(v)\leq b$ for each $v\in V(G)$. A $[k,k]$-factor is a regular factor, which is usually called a \textit{$k$-factor}. The research on the existence of various factors can be dated back to Hall's marriage theorem~\cite{Hall35}, Tutte's $1$-factor theorem~\cite{Tutte47} and Tutte's $k$-factor theorem~\cite{Tutte52}. For an integer-valued function $f$ defined on $V(G)$, an \textit{$f$-factor} of a graph $G$ is a spanning subgraph $H$ of $G$ such that $d_{H}(v)=f(v)$ for all $v\in V(G)$. In fact, Tutte~\cite{Tutte52} provided a necessary and sufficient condition for the existence of an $f$-factor. Ore~\cite{Ore57} then studied $f$-factors in bipartite graphs, which implies a necessary and sufficient condition for the existence of a biregular bipartite spanning subgraphs (see Theorem~\ref{thm::ore} and Corollary~\ref{cor::biregular} in Section~\ref{sec::bifactor} for the results).
%as well as Lov\'asz's general factor theorem~\cite{Lov70}.

An $(a,b)$-biregular bipartite spanning subgraph of a bipartite graph $G$ is also called an \textit{$(a,b)$-biregular factor} of $G$. We intend to investigate biregular factors in bipartite graphs via eigenvalues. For a graph $G$ on $n$ vertices, let $A(G)$ denote  the adjacency matrix of $G$ and let $\lambda_i:=\lambda_i(G)$ denote the $i$-th largest eigenvalue of $A(G)$ for $i=1,2,\ldots, n$. Specifically, the largest eigenvalue of $A(G)$ is called the \textit{spectral radius} of $G$ and is denoted by $\rho(G)$.

The existences of various factors have been well studied from eigenvalues. In 2005, Brouwer and Haemers \cite{A.B} gave a sufficient condition for a regular graph to contain a 1-factor in terms of the third largest adjacency eigenvalue. Their result was improved in \cite{S.C-1,S.C,S.C-2} and extended to $k$-factors in \cite{H.L-3,H.L-4, Gu14}. O \cite{S.O} considered the graphs that are not necessarily regular, and  provided a spectral condition to guarantee the existence of a 1-factor in a graph. Subsequently, Cho, Hyun, O and Park \cite{E.C} proposed a conjecture on the existence of $[a,b]$-factors in a graph of order $n$ via its spectral radius. This conjecture was completely confirmed in \cite{WZ2023}.

An $(X, Y)$-bipartite graph is called \emph{balanced} if $|X| = |Y|$. Fan and Lin \cite{FL22} gave a sufficient condition for the existence of a $k$-factor in a balanced bipartite graph via spectral radius.  In this paper, we will extend the result in \cite{FL22} by exploring the existence of an $(a,b)$-biregular factor in a bipartite graph.

Denote by $G\backslash E(H)$ the graph obtained from $G$ by deleting the edges of $H$, where $H$ is a subgraph of $G$.
\begin{theorem}\label{thm::biregular-factor}
Suppose that $a$ and $b$ are two positive integers with $b>a$. Let $G$ be a connected $(X,Y)$-bipartite graph with $am=bn$, where $|X|=m$ and $|Y|=n$. If $m\geq n+b$ and
$$\rho(G)\geq\rho(K_{m,n}\backslash E(K_{1,n-a+1})),$$
then $G$ contains an $(a,b)$-biregular factor, unless $G\cong K_{m,n}\backslash E(K_{1,n-a+1})$.
\end{theorem}

%A basic result in graph theory asserts that every connected graph has a spanning tree. 
We also study various spanning trees in bipartite graphs. In the past several decades, many researchers have explored the existence of spanning trees under some specified conditions. A \textit{$k$-tree} is a spanning tree with every vertex of degree at most $k$, where $k\geq 2$ is an integer. In 2010, Fiedler and Nikiforov \cite{Fiedler-Nikiforov} gave spectral radius conditions for the existence of a 2-tree in a graph, which was improved in \cite{Ning-Ge,Li-Ning}. For $k\geq 3$, Wong \cite{Wong2013} showed that for any $d$-regular connected graph $G$, if $\lambda_4(G) <d-\frac{d}{(k-2)(d+1)}$, then $G$ has a $k$-tree. Later, this result was generalized by  Cioab\u{a} and Gu \cite{Cioaba-Gu}, who incorporated the edge connectivity into their analysis. Recently, a spectral radius condition was obtained in \cite{FGHL22}. Up to now, much attention has been paid to the relationship between the eigenvalues and the spanning tree under some specified conditions of a graph, we refer the reader to  \cite{ALYL22,ALY23,MALWH23}. For bipartite graphs, Frank and Gy\'arf\'as~\cite{FG76}, Kaneko and Yoshimoto~\cite{KY00} independently investigated the condition for the existence of a spanning tree $T$ in a connected $(X,Y)$-bipartite graph such that $d_T(u)\ge k$ for all $u\in X$.  In this paper, we study the existence of such spanning trees by means of the spectral radius.

\begin{theorem}\label{thm::spanning-tree1}
Suppose that $k\geq 3$ is an integer, and $G$ is a connected $(X,Y)$-bipartite graph with $|X|=m$ and $|Y|=n$. If $n>(k-1)m$, $m\geq k+1$ and
$$\rho(G)\geq \rho(K_{m,n}\backslash E(K_{1,n-k+1})),$$
 then $G$ has a spanning tree $T$ such that $d_T(u)\ge k$ for all $u\in X$, unless $G\cong K_{m,n}\backslash E(K_{1,n-k+1})$.
\end{theorem}
\begin{remark}
For a connected $(X,Y)$-bipartite graph $G$, if $n\leq (k-1)m$, let $S=X$, we can deduce that $G$ contains no spanning trees $T$ such that $d_T(u)\ge k$ for all $u\in X$ by Lemma \ref{lem::tree-condition}. Thus, we need consider $n>(k-1)m$ in Theorem \ref{thm::spanning-tree1}. 
\end{remark}

%By letting $k=2$ in Theorem~\ref{thm::spanning-tree1}, we obtain a spectral radius condition for the existence of a spanning tree whose leaves are all contained in $Y$ in a connected $(X,Y)$-bipartite graph. 
For $k=2$ in the above mentioned spanning tree, it means all leaves are contained in $Y$. For this special type of spanning tree, we can even do better by imposing the minimum degree as a parameter, as below.

\begin{theorem}\label{thm::spanning-tree2}
Let $G$ be a connected $(X,Y)$-bipartite graph with minimum degree $\delta$, where $|X|=m$ and $|Y|=n$. If $n>m\geq 3$ and
$$\rho(G)\geq \rho(K_{m,n}\backslash E(K_{\delta,n-\delta})),$$
 then $G$ has a spanning tree whose leaves are all contained in $Y$, unless $G\cong K_{m,n}\backslash E(K_{\delta,n-\delta})$.
\end{theorem}

The \textit{spanning tree packing number} of a graph $G$, denoted by $\tau(G)$, is the maximum number of edge-disjoint spanning trees contained in $G$.  
Cioab\u{a}~\cite{Cioa10} proved that for any $d$-regular connected graph $G$, if $\lambda_{2}(G)<d-\frac{2(2k-1)}{d+1}$ for $d\geq 2k\geq 4$, then $\tau(G)\geq k$. Subsequently, Cioab\u{a} and Wong \cite{CiWo12} further conjectured that the sufficient condition can be improved to $\lambda_{2}(G)<d-\frac{2k-1}{d+1}$, and they did the preliminary work of this conjecture for $k=2,3$. Gu et al.~\cite{GLLY12} extended it to general graphs, and obtained a partial result. The conjecture  in \cite{GLLY12} was completely settled in \cite{LHGL14}, which also implied the truth of the conjecture of Cioab\u{a} and Wong \cite{CiWo12}. Recently, the results of \cite{LHGL14} have been shown to be essentially best possible in \cite{coppww22} by constructing extremal graphs. The spanning tree packing number was studied via the spectral radius in \cite{FGL23}, in which we proved that for a connected graph $G$ with minimum degree $\delta\geq 2k\geq 4$ and order $n\geq 2\delta+3$, if $\rho(G)\geq \rho(B_{n,\delta+1}^{k-1})$, then $\tau(G)\geq k$ unless $G\cong B_{n,\delta+1}^{k-1}$, where $B_{n,\delta+1}^{k-1}$ is a graph obtained from $K_{\delta+1}\cup K_{n-\delta-1}$ by adding $k-1$ edges joining a vertex in $K_{\delta+1}$ and $k-1$ vertices in $K_{n-\delta-1}$. This motivated us to study an analogue in bipartite graphs. We prove the following result.

\begin{theorem}\label{thm::edge-spanningtree}
Let $n$ and $k$ be two positive integers with $n\geq 2k+2$, and let $G$ be a connected balanced bipartite graph of order $2n$. If 
$$\rho(G)\geq \rho(K_{n,n}\backslash E(K_{1,n-k+1})),$$
then $\tau(G)\geq k$, unless $G\cong K_{n,n}\backslash E(K_{1,n-k+1})$.
\end{theorem}

\section{Preliminaries}
Let $e(G)$ denote the number of edges in $G$.
\begin{lemma}
[See \cite{Nosal}]\label{lem::edge}
Let $G$ be a bipartite graph. Then
$$\rho(G)\leq \sqrt{e(G)},$$
with equality if and only if $G$ is a complete bipartite graph with possibly some isolated vertices.
\end{lemma}

Let $M$ be a real $n\times n$ matrix, and let $X=\{1,2,\ldots,n\}$. Given a partition $\Pi:X=X_{1}\cup X_{2}\cup \cdots \cup X_{k}$,  the matrix $M$ can be correspondingly partitioned as
$$
M=\left(\begin{array}{ccccccc}
M_{1,1}&M_{1,2}&\cdots &M_{1,k}\\
M_{2,1}&M_{2,2}&\cdots &M_{2,k}\\
\vdots& \vdots& \ddots& \vdots\\
M_{k,1}&M_{k,2}&\cdots &M_{k,k}\\
\end{array}\right).
$$
The \textit{quotient matrix} of $M$ with respect to $\Pi$ is defined as the $k\times k$ matrix $B_\Pi=(b_{i,j})_{i,j=1}^k$ where $b_{i,j}$ is the  average value of all row sums of $M_{i,j}$.
The partition $\Pi$ is called \textit{equitable} if each block $M_{i,j}$ of $M$ has constant row sum $b_{i,j}$.
Also, we say that the quotient matrix $B_\Pi$ is \textit{equitable} if $\Pi$ is an equitable partition of $M$.

\begin{lemma}[Brouwer and Haemers \cite{BH}, p. 30; Godsil and Royle \cite{C.Godsil}, pp. 196--198]\label{lem::equ}
Let $M$ be a real symmetric matrix, and let $\lambda_{1}(M)$ be the largest eigenvalue of $M$. If $B_\Pi$ is an equitable quotient matrix of $M$, then the eigenvalues of  $B_\Pi$ are also eigenvalues of $M$. Furthermore, if $M$ is nonnegative and irreducible, then $\lambda_{1}(M) = \lambda_{1}(B_\Pi).$
\end{lemma}

\section{The proofs of Theorems \ref{thm::biregular-factor}, \ref{thm::spanning-tree1} and \ref{thm::spanning-tree2}}\label{sec::bifactor}

\begin{lemma} \label{ineuq}
 Let $x,y$ and $a$ be three positive integers with $x>y\geq a$. Then
 $$\rho(K_{x,y}\backslash E(K_{1,y-a+1}))>\rho(K_{x,y}\backslash E(K_{x-a+1,1})).$$
\end{lemma}
\renewcommand\proofname{\bf Proof}
\begin{proof}
Note that $A(K_{x,y}\backslash E(K_{1,y-a+1}))$ has the equitable quotient matrix
$$
A=\begin{bmatrix}
0 &0 & a-1& 0\\
0&0&a-1 &y-a+1\\
1&x-1& 0& 0\\
0&x-1&0 &0
\end{bmatrix}.
$$
By a simple computation, the characteristic polynomial of $A$ is
\begin{equation*}
\begin{aligned}
\varphi(\lambda)=\lambda^4 + (-xy-a+y+1)\lambda^2-a^2x+axy+a^2+2xa-ay-xy-2a - x + y + 1.
\end{aligned}
\end{equation*}
Also, note that $A(K_{x,y}\backslash E(K_{x-a+1,1}))$ has the equitable quotient matrix
$$
B=\begin{bmatrix}
0 &0 & 1& y-1\\
0&0&0 &y-1\\
a-1&0& 0& 0\\
a-1&x-a+1&0 &0
\end{bmatrix}.
$$
By a simple computation, the characteristic polynomial of $B$ is
\begin{equation*}
\begin{aligned}
\phi(\lambda)=\lambda^4 + (-xy - a + x + 1)\lambda^2 - a^2y + axy + a^2 - xa + 2ay - xy - 2a + x - y + 1.
\end{aligned}
\end{equation*}
For $\lambda\geq \sqrt{(x-1)y}$, we get
$$\phi(\lambda)-\varphi(\lambda)=(x-y)(a^2\!+\!\lambda^2\!-\!3a\!+\!2)\geq (x\!-\!y)(a^2\!+\!(x\!-\!1)y\!-\!3a\!+\! 2)\geq (x\!-\!y)((2a\!-\!1)(a\!-\!1)\!+\!1)>0$$
due to $x>y\geq a+1$ and $a\geq 1$, which leads to
$\varphi(\lambda)<\phi(\lambda)$. Since $K_{x-1,y}$ is a proper spanning subgraph of  $K_{x,y}\backslash E(K_{1,y-a+1})$, we have $\rho(K_{x,y}\backslash E(K_{1,y-a+1}))>\rho(K_{x-1,y})=\sqrt{(x-1)y}$. It follows that $\lambda_1(K_{x,y}\backslash E(K_{1,y-a+1}))>\lambda_1(K_{x,y}\backslash E(K_{x-a+1,1}))$. Combining this with Lemma \ref{lem::equ}, we get
$$\rho(K_{x,y}\backslash E(K_{1,y-a+1}))>\rho(K_{x,y}\backslash E(K_{x-a+1,1})).$$ 
This completes the proof.
\end{proof}
For $X,Y\subset V(G)$, we denote  by $e_{G}(X,Y)$ the number of edges with one endpoint in $X$ and one endpoint in $Y$.
\begin{theorem}[Ore~\cite{Ore57}] %Folkman and Fulkerson~\cite{FF70}
\label{thm::ore}
Let $G$ be an $(X, Y)$-bipartite graph and let $f:V(G)\rightarrow \{1,2,3,4,\ldots\}$ be a function. Then $G$ has an $f$-factor if and only if $\sum_{u\in X}f(u) =\sum_{u\in Y}f(u)$ and
\begin{eqnarray*}
\delta(S,T) 
&=& \sum_{u\in T}f(u) + \sum_{u\in S}\left(d_{G-T}(u) - f(u)\right)\\
&=& \sum_{u\in T}f(u) + \sum_{u\in S}\left(d_{G}(u) - f(u)\right) -e_{G}(S,T)\\
&=& \sum_{u\in T}f(u) - \sum_{u\in S}f(u) + e_{G}(S,Y-T)\ge 0
\end{eqnarray*}
for all subsets $S\subseteq X$ and $T\subseteq Y$.
\end{theorem}
%Note that in the above theorem, instead of using $\delta(S,T)$, it is equivalent to use
%\begin{eqnarray*}
%\delta(T,S) 
%&=& \sum_{u\in S}f(u) + \sum_{u\in T}\left(d_{G-S}(u) - f(u)\right)\\
%&=& \sum_{u\in S}f(u) + \sum_{u\in T}\left(d_{G}(u) - f(u)\right) -e(S,T)\\
%&=& \sum_{u\in S}f(u) - \sum_{u\in T}f(u) + e(X-S,T)\ge 0.
%\end{eqnarray*}
\begin{corollary}\label{cor::biregular}
An $(X,Y)$-bipartite graph $G$ has an $(a,b)$-biregular factor if and only if $a|X|=b|Y|$ and $$\delta(S,T)=e_{G}(S,Y-T) +b|T| -a|S|\ge 0$$ for all subsets $S\subseteq X$ and $T\subseteq Y$.
\end{corollary}
\begin{lemma}[Frank and Gy\'arf\'as~\cite{FG76}, Kaneko and Yoshimoto~\cite{KY00}]\label{lem::tree-condition}
Let $G$ be a connected $(X,Y)$-bipartite simple graph and $f:X\rightarrow \{2,3,4,\ldots\}$ be a function. Then $G$ has a spanning tree $T$ such that $d_T(u)\ge f(u)$ for all $u\in X$ if and only if 
$$|N_G(S)|\ge \sum_{u\in S}f(u) -|S| +1$$
for all nonempty subsets $S\subseteq X$. In particular, if $|N_G(S)|\ge |S|+1$ for all nonempty subsets $S\subseteq X$, then $G$ has a spanning tree whose leaves are all contained in $Y$.
\end{lemma}

\begin{lemma}\label{lem::delta}
Let $x$, $y$, $\delta$ and $s$ be positive integers with $y>x\geq s+\delta$ and $s\geq \delta+1$. Then
$$\rho(K_{x,y}\backslash E(K_{\delta,y-\delta}))> \rho(K_{x,y}\backslash E(K_{s,y-s})).$$
\end{lemma}

\renewcommand\proofname{\bf Proof}
\begin{proof}
 Observe that $K_{x,y}\backslash E(K_{s,y-s})$ has the equitable quotient matrix
$$
B_\Pi^s=\begin{bmatrix}
0 &0 & s& 0\\
0 &0 &s&y-s \\
  s &x-s &0&0\\
  0 &x-s &0&0\\ 
\end{bmatrix}.
$$
By a simple calculation, the characteristic polynomial of $B_\Pi^s$ is
\begin{equation*}
\begin{aligned}
\varphi(B_{\Pi}^s,\lambda)=\lambda^4+(-s^2+sy-xy)\lambda^2+s^4-s^3x-s^3y+s^2xy.
\end{aligned}
\end{equation*}
Also, note that $A(K_{x,y}\backslash E(K_{\delta,y-\delta})$ has the equitable quotient matrix $B_{\Pi}^{\delta}$, which is obtained by replacing $s$ with $\delta$ in $B_\Pi^s$. Then
\begin{eqnarray*}
\varphi(B_{\Pi}^s,\lambda)\!-\!\varphi(B_{\Pi}^{\delta},\lambda)=(s\!-\!\delta)((y-\!\delta\!-\!s)\lambda^2\!+\!\delta^2(\delta\!+\!s\!-\!x\!-\!y)\!+\!\delta(s^2\!-\!sx\!-\! sy\!+\! xy)\!+\!s^2(s\!-\!x\!-\!y)\!+\!sxy).
\end{eqnarray*}
Let $g(\lambda)=(y-\!\delta\!-\!s)\lambda^2\!+\!\delta^2(\delta\!+\!s\!-\!x\!-\!y)\!+\!\delta(s^2\!-\!sx\!-\! sy\!+\! xy)\!+\!s^3\!-\!s^2x\!-\!s^2y\!+\!sxy$. For $\lambda\geq \sqrt{(x-\delta)y}$, we get
\begin{equation*}
\begin{aligned}
  g(\lambda)&=(y-\!\delta\!-\!s)\lambda^2\!+\!\delta^2(\delta\!+\!s\!-\!x\!-\!y)\!+\!\delta(s^2\!-\!sx\!-\! sy\!+\! xy)\!+\!s^3\!-\!s^2x\!-\!s^2y\!+\!sxy\\
    &\geq (y^2-s^2\!-\!\delta^2\!-\! s\delta )x\!+\!\delta^3\!+\!s\delta^2\!+\!(s^2\!-\!y^2)\delta\!+\!(s\!-\!y)s^2~~(\mbox{since $\lambda\geq \sqrt{(x\!-\!\delta)y}$ and $y> s\!+\!\delta$})\\
     &\geq (y^2-s^2\!-\!\delta^2\!-\! s\delta )(\delta+s)\!+\!\delta^3\!+\!s\delta^2\!+\!(s^2\!-\!y^2)\delta\!+\!(s\!-\!y)s^2~~(\mbox{since $x\geq s\!+\!\delta$})\\
     &= s(\delta + y)(y-\delta- s)\\
     &> 0 ~~~(\mbox{since $y> s+\delta$}).
\end{aligned}
\end{equation*}
Combining this with $s\geq \delta+1$, we deduce that
$\varphi(B_{\Pi}^s,\lambda)>\varphi(B_{\Pi}^{\delta},\lambda)$. Since $K_{x,y}\backslash E(K_{\delta,y-\delta})$ contains $K_{x-\delta,y}$ as a proper subgraph, we have $\rho(K_{x,y}\backslash E(K_{\delta,y-\delta}))>\sqrt{(x-\delta)y}$ , and hence $\lambda_1(B_{\Pi}^{\delta})>\lambda_1(B_{\Pi}^{s})$. Combining this with Lemma \ref{lem::equ}, we can deduce that $\rho(K_{x,y}\backslash E(K_{\delta,y-\delta}))>\rho(K_{x,y}\backslash E(K_{s,y-s}))$.
\end{proof}

Now we give the proofs of Theorems \ref{thm::biregular-factor}, \ref{thm::spanning-tree1} and \ref{thm::spanning-tree2}.
\renewcommand\proofname{\bf Proof of Theorem \ref{thm::biregular-factor}}
\begin{proof}
Suppose to the contrary that $G$ contains no $(a,b)$-biregular factors. By Corollary \ref{cor::biregular}, there exist subsets $S\subseteq X$ and $T\subseteq Y$ such that $e_{G}(S,Y-T)+b|T|-a|S|\leq -1$. Let $|S|=s$ and $|T|=t$. Then
\begin{equation}\label{equ::biregular-1}
e_{G}(S,Y-T)\leq as-bt-1.
\end{equation}
We have the following three claims.

{\flushleft\bf{Claim 1.}} $s\geq t+1.$

Otherwise, $s\leq t$. Combining this with $b>a$ and (\ref{equ::biregular-1}), we have
$$0\leq e_{G}(S,Y-T)\leq as-bt-1<bt-bt-1=-1,$$
a contradiction. Hence, $s\geq t+1$.\qed

{\flushleft\bf{Claim 2.}} $d_{G}(v)\geq a$ for all $v\in V(G)$.

Otherwise, there exists a vertex $w\in V(G)$ such that $d_{G}(w)\leq a-1$. If $w\in X$, then $G$ is a spanning subgraph of $K_{m,n}\backslash E(K_{1,n-a+1})$. It follows that 
$$\rho(G)\leq \rho(K_{m,n}\backslash E(K_{1,n-a+1})),$$ with equality if and only if $G\cong K_{m,n}\backslash E(K_{1,n-a+1})$, a contradiction. If $w\in Y$, then $G$ is a spanning subgraph of $K_{m,n}\backslash E(K_{m-a+1,1})$. Combining this with Lemma \ref{ineuq}, we have 
$$\rho(G)\leq \rho(K_{m,n}\backslash E(K_{m-a+1,1}))<\rho(K_{m,n}\backslash E(K_{1,n-a+1})),$$
which also leads to a contradiction. This implies that $d_{G}(v)\geq a$ for all $v\in V(G)$. \qed

{\flushleft\bf{Claim 3.}} $1\leq t\leq n-1$.

Otherwise, $T=Y$ or $T=\emptyset$. If $T=Y$, then $t=n$. Combining this with  (\ref{equ::biregular-1}), $s\leq m$ and $am=bn$, we get 
$$0=e_{G}(S,Y-T)=e_{G}(S,\emptyset)\leq as-bn-1=a(s-m)-1\le-1,$$
a contradiction. If $T=\emptyset$, then $t=0$. By Claim 2 and (\ref{equ::biregular-1}), we have 
$$as\leq e_{G}(S,Y-T)=e_{G}(S,Y)\leq as-1,$$
which also leads to a contradiction. Thus, $1\leq t\leq n-1$. \qed

Observe that $d_{Y-T}(v)\geq a-t$ for each $v\in S$. Then again by (\ref{equ::biregular-1}), we have 
$$(a-t)s\leq e_{G}(S,Y-T)\leq as-bt-1,$$
and hence $s\geq b+\frac{1}{t}$. Furthermore, we can deduce that $s\geq b+1$ because $s$ is a positive integer and $t\geq 1$. On the one hand, note that $\rho(G)\geq \rho(K_{m,n}\backslash E(K_{1,n-a+1}))>\rho(K_{m-1,n})=\sqrt{(m-1)n}$. Then by Lemma \ref{lem::edge}, we have
\begin{equation}\label{equ::biregular-2}
e(G)>(m-1)n.
\end{equation}
On the other hand, by  (\ref{equ::biregular-1}), we can deduce that
\begin{equation}\label{equ::biregular-3}
\begin{aligned}
  e(G)&= e_{G}(S,T)+e_{G}(S,Y-T)+e_{G}(X-S,Y)\\
      &\leq st+as-bt-1+(m-s)n ~~(\mbox{by (\ref{equ::biregular-1})})\\
      &=(m-1)n-((s-1)n-(s-b)t-as+1).
\end{aligned}
\end{equation}
Let $h(s)=(s-1)n-(s-b)t-as+1$. Then we consider the following two cases.

{\flushleft\bf{Case 1.}} $s\leq n$.

Note that $s\geq b+1$. Then by Claim 1 and $s\leq n$, we get
\begin{equation*}
\begin{aligned}
  h(s)&=(s-1)n-(s-b)t-as+1\\
      &\geq (n-s+b-a-1)(s-1)+s-a~~(\mbox{since $t\leq s-1$ and $s\geq b+1$})\\
      &>0 ~~(\mbox{since $s\leq n$, $b\geq a+1$ and $s\geq b+1$}).\\
     \end{aligned}
\end{equation*}
Combining this with (\ref{equ::biregular-3}), we deduce that $e(G)<(m-1)n$, which contradicts (\ref{equ::biregular-2}).

{\flushleft\bf{Case 2.}} $s> n$.

If $n\geq t+a+1$, then 
$$h(s)=s(n-t-a)-n+bt+1>s(n-t-a-1)+bt+1>0$$
due to $s>n$. If $n\leq t+a$, then \begin{equation*}
\begin{aligned}
  h(s)&=s(n-t-a)-n+bt+1\\
      &\geq m(n-t-a)-n+bt+1 ~~(\mbox{since $s\leq m$ and $n\leq t+a$})\\
       &=(m-b)(n-t)-n+1 ~~(\mbox{since $am=bn$})\\
       &\geq m-b-n+1 ~~(\mbox{since $t\leq n-1$})\\
       &>0  ~~(\mbox{since $m\geq b+n$}).
     \end{aligned}
\end{equation*}    
Thus, we can deduce that $h(s)>0$ for $s>n$. Again by (\ref{equ::biregular-3}), we have $e(G)<(m-1)n$, which also contradicts (\ref{equ::biregular-2}). 

This completes the proof.
\end{proof}

\renewcommand\proofname{\bf Proof of Theorem \ref{thm::spanning-tree1}}
\begin{proof}
Suppose to the contrary that $G$ has no spanning trees $T$ such that $d_T(u)\ge k$ for all $u\in X$. Then by Lemma \ref{lem::tree-condition}, there exists some nonempty subset $S\subseteq X$ with $|S|=s$ such that $|N_G(S)|\leq (k-1)s$.  
If $s=1$, then $G$ is a spanning subgraph of $K_{m,n}\backslash E(K_{1,n-k+1})$. Hence, $\rho(G)\leq \rho(K_{m,n}\backslash E(K_{1,n-k+1}))$, with equality if and only if $G\cong K_{m,n}\backslash E(K_{1,n-k+1})$, a contradiction. Next, we consider $s\geq 2$.
Since $\rho(G)\geq \rho(K_{m,n}\backslash E(K_{1,n-k+1}))$ and $K_{m-1,n}$ is a proper subgraph of $K_{m,n}\backslash E(K_{1,n-k+1})$, we have $$\rho(G)\geq \rho(K_{m,n}\backslash E(K_{1,n-k+1}))>\sqrt{(m-1)n}.$$
Combining this with Lemma \ref{lem::edge}, we get 
\begin{equation}\label{equ::tree-1}
e(G)>(m-1)n.
\end{equation}
Note that $G$ is a spanning subgraph of $K_{m,n}\backslash E(K_{s,n-(k-1)s})$. Then 
\begin{equation}\label{equ::tree-2}
\begin{aligned}
  e(G)&= e(S,N_{G}(S))+e(X-S,Y)\\
     &\leq (k-1)s^2+(m-s)n\\
     &=(m-1)n-((s-1)n-(k-1)s^2).
\end{aligned}
\end{equation}
Let $f(s)=(s-1)n-(k-1)s^2$. If $2\leq s\leq k$, by a simple computation, we get $f(2)=n-4(k-1)\geq (k-2)^2> 0$ and  $f(k)=(k-1)(n-k^2)\geq 0$ due to $n\geq (k-1)m+1$, $m\geq k+1$ and $k\geq 3$. Thus, $f(s)\geq \min\{f(2),f(k)\}\geq 0$ for $2\leq s\leq k$. Combining this with (\ref{equ::tree-2}), we get $e(G)\leq (m-1)n-f(s)\leq (m-1)n$, which contradicts (\ref{equ::tree-1}). Thus, we consider $s\geq k+1$ in the following. We assert that $S$ is a proper subset of $X$. Otherwise, $S=X$. Then $|N_{G}(X)|\leq (k-1)m$. It follows that $n-|N_{G}(X)|\geq (k-1)m+1-(k-1)m=1$. This is impossible because $G$ is connected. Thus, $m\geq s+1$, and hence $n\geq (k-1)m+1\geq (k-1)s+k$. By (\ref{equ::tree-2}), we have 
\begin{equation*}
\begin{aligned}
  e(G)&\leq (m-1)n-((s-1)n-(k-1)s^2)\\
     &\leq (m-1)n-((s-1)((k-1)s+k)-(k-1)s^2)~~~(\mbox{since $n\geq (k-1)s+k$ and $s\geq k+1\geq 4$})\\
     &=(m-1)n-(s-k)\\
     &<(m-1)n ~~~(\mbox{since $s\geq k+1$}),
\end{aligned}
\end{equation*}
a contradiction. 

This completes the proof.
\end{proof}
%\begin{problem}
%Study the spectral radius and bounded-degree spanning trees of bipartite graphs.
%\end{problem}

\renewcommand\proofname{\bf Proof}

\renewcommand\proofname{\bf Proof of Theorem \ref{thm::spanning-tree2}}
\begin{proof}
Suppose to the contrary that $G$ contains no spanning trees whose leaves are all contained in $Y$. Then by Lemma \ref{lem::tree-condition}, there exists some nonempty subset $S\subseteq X$ with $|S|=s$ such that $|N_G(S)|\leq s$. Observe that $G$ is a spanning subgraph of $K_{m,n}\backslash E(K_{s,n-s})$. Thus, 
\begin{equation}\label{equ::tree2-1}
\begin{aligned}
\rho(G)\leq \rho(K_{m,n}\backslash E(K_{s,n-s})),
\end{aligned}
\end{equation}
with equality if and only if $G\cong K_{m,n}\backslash E(K_{s,n-s})$. 
It is easy to verify that $s\geq \delta$. If $s\geq\delta+1$, by $n>m\geq \delta+s$ and Lemma \ref{lem::delta}, we have
\begin{equation*}
\begin{aligned}
\rho(K_{m,n}\backslash E(K_{s,n-s}))<\rho(K_{m,n}\backslash E(K_{\delta,n-\delta})),
\end{aligned}
\end{equation*}
a contradiction.
If $s=\delta$, by (\ref{equ::tree2-1}), we have 
$$\rho(G)\leq \rho(K_{m,n}\backslash E(K_{\delta,n-\delta})),$$
with equality if and only if $G\cong K_{m,n}\backslash E(K_{\delta,n-\delta})$, which also leads to a contradiction.

This completes the proof.
\end{proof}

\section{The proof of Theorem \ref{thm::edge-spanningtree}}

\begin{lemma}\label{lem::1}
For $1\leq i\leq s$, let $a_i$ and $b_i$ be non-negative integers with  $a_i+b_i\geq 2$, $\sum_{1\leq i\leq s}a_i=a$ and $\sum_{1\leq i\leq s}b_i=b$. If $b\geq a$, then 
$$\sum_{1\leq i\leq s}a_ib_i\leq a(b-(s-1)).$$
\end{lemma}
\renewcommand\proofname{\bf Proof}
\begin{proof}
Suppose that $p=|\{j|~a_j=0, 1\leq j\leq s\}|$ and $q=|\{j|~b_j=0, 1\leq j\leq s\}|$. 
If $p=q=0$, then $a_i\geq 1$ and $b_i\geq 1$ for all $1\leq i\leq s$, and hence
\begin{equation*}
\begin{aligned}
 a(b-(s-1))-\sum_{1\leq i\leq s}a_ib_i&= \Big(\sum_{1\leq i\leq s}a_i\Big)\Big(\sum_{1\leq i\leq s}b_i-(s-1)\Big)-\sum_{1\leq i\leq s}a_ib_i\\
     &=\sum_{1\leq j\leq s}\sum_{1\leq i\neq j\leq s}a_j(b_i-1)\\
     &\geq 0,
\end{aligned}
\end{equation*}
as required. 
Next, we consider $\max\{p,q\}\geq 1$. Without loss of generality, we may assume that $a_i=0$ and $b_j=0$ for $1\leq i\leq p$ and $p+1\leq j\leq p+q$. Then
$b_i\geq 2$ and $a_j\geq 2$ for $1\leq i\leq p$ and $p+1\leq j\leq p+q$. Using the results above, we get
\begin{equation*}
\begin{aligned}
\sum_{1\leq i\leq s}a_ib_i&=\sum_{p+q+1\leq i\leq s}a_ib_i\\
&\leq \max\Big\{\Big(a\!-\!\sum_{p\!+\!1\leq i\leq p\!+\!q}a_i\Big)\Big(b\!-\!\sum_{1\leq i\leq p}b_i\!-\!(s\!-\!p\!-\!1)\Big), \Big(b\!-\!\sum_{1\leq i\leq p}b_i\Big)\Big(a\!-\!\sum_{p\!+\!1\leq i\leq p\!+\!q}a_i\!-\!(s\!-\!q\!-\!1)\Big)\Big\}\\
&\leq \max\{(a-2q)(b-p-(s-1)), (b-2p)(a-q-(s-1))\}\\
&~~(\mbox{since $b_i\geq 2$ and $a_j\geq 2$ for $1\leq i\leq p$ and $p+1\leq j\leq p+q$})\\
&<a(b-(s-1))~~(\mbox{since $\max\{p,q\}\geq 1$ and $b\geq a$}).
\end{aligned}
\end{equation*}
Thus the result follows.
\end{proof}

Let $(V_1,V_2,\cdots,V_t)$ be a partition of $V(G)$, and let $e_{G}(V_1,V_2,\cdots,V_t)$  denote the number of edges in $G$ whose endpoints lie in different parts of $(V_1,V_2,\cdots,V_t)$. A part is \textit{trivial} if it contains a single vertex. The following fundamental theorem on spanning tree packing number of a graph was established by Nash-Williams \cite{Nash-Williams} and Tutte \cite{Tutte}, independently.
\begin{theorem}[Nash-williams~\cite{Nash-Williams} and Tutte~\cite{Tutte}]
\label{lem::Nash-williams-Tutte}
Let $G$ be a connected graph. Then $\tau(G)\geq k$ if and only if for any partition $(V_1,V_2,\cdots,V_t)$ of $V(G)$, $$e_{G}(V_1,V_2,\cdots, V_t)\geq k(t-1).$$
\end{theorem}

For $X\subseteq V(G)$, let $G[X]$ be the subgraph of $G$ induced by $X$, and let $e(X)$ be the number of edges in $G[X]$.
\renewcommand\proofname{\bf Proof of Theorem \ref{thm::edge-spanningtree}}
\begin{proof}
Let $G=(X,Y)$ be a balanced bipartite graph of order $2n$ with $\tau(G)\leq k-1$. By Lemma \ref{lem::Nash-williams-Tutte}, there exists a partition $(V_1,V_2,\cdots, V_t)$ of $V(G)$ such that
\begin{equation}\label{equ::edgetree-1}
\begin{aligned}
e_{G}(V_1,V_2,\cdots, V_t)\leq k(t-1)-1.
\end{aligned}
\end{equation}
Since $\rho(G)\geq \rho(K_{n,n}\backslash E(K_{1,n-k+1}))$ and $K_{n,n-1}$ is a proper subgraph of $K_{n,n}\backslash E(K_{1,n-k+1})$, we have 
$$\rho(G)\geq \rho(K_{n,n}\backslash E(K_{1,n-k+1}))>\sqrt{n(n-1)}.$$ Combining this with Lemma \ref{lem::edge}, we have
\begin{equation}\label{equ::edgetree-2}
e(G)>n(n-1).
\end{equation}
We assume that $|V_1|\geq |V_2|\geq \cdots\geq |V_t|\geq 1$. Let $|V_i\cap X|=a_i$ and $|V_i\cap Y|=b_i$ where $1\leq i\leq s$. We divide the proof into the following two cases.

{\flushleft\bf Case 1.} $t=2$.

In this case, we assert that $|V_2|=1$. Otherwise, $|V_2|\geq 2$. Since $|V_1|\geq |V_2|$, it follows that $a_1+b_1\geq n$. Without loss of generality, we may assume that $a_1\geq b_1$. Take $a_1=b_1+p$ where $p\geq 0$. Thus, $\frac{n+p}{2}\leq a_1\leq n$. If $a_1=n$, then $a_2=n-a_1=0$, and hence $V_2\subseteq Y$. Observe that $G$ is connected and $\sum_{u\in V_2}d_{G}(u)=e_{G}(V_1,V_2)\leq k-1$ due to (\ref{equ::edgetree-1}). Thus, $d_{G}(u)<k-1$ for all $u\in V_2$. This implies that $G$ is a proper spanning subgraph of $K_{n,n}\backslash E(K_{1,n-k+1})$. Thus, 
$$\rho(G)< \rho(K_{n,n}\backslash E(K_{1,n-k+1})),$$ 
a contradiction. Thus, $\frac{n+p}{2}\leq a_1\leq n-1$. Combining this with (\ref{equ::edgetree-1}), we get
\begin{equation*}
\begin{aligned}
e(G)&=e(V_1)+e(V_2)+e_{G}(V_1,V_2)\\
&\leq a_1b_1+a_2b_2+k-1\\
     &=a_1b_1+(n-a_1)(n-b_1)+k-1~~(\mbox{since $a_2=n-a_1$ and $b_2=n-b_1$})\\
     &=n(n-1)-((a_1+b_1-1)n-2a_1b_1-k+1)\\
     &=n(n-1)-((2b_1\!+\!p\!-\!1)(a_1\!+\!a_2)\!-\!2a_1b_1\!-\!k+1)~~(\mbox{since $n=a_1+a_2$ and $a_1=b_1+p$})\\
     &=n(n-1)-(2a_2b_1+(p-1)n-k+1)\\
     &=n(n-1)-(2(n-a_1)(a_1-p)+(p-1)n-k+1)\\
     &=n(n-1)-(2(-a^2_1+(n+p)a_1-pn)+(p-1)n-k+1)\\
     &\leq n(n-1)-(2(-(n\!-\!1)^2\!+\!(n\!+\!p)(n\!-\!1)\!-\!pn)\!+\!(p\!-\!1)n\!-\!k\!+\!1)~~(\mbox{since $\frac{n+p}{2}\leq a_1\leq n-1$})\\
     &=n(n-1)-((n-2)p+n-k-1)\\
     &< n(n-1) ~~(\mbox{since $n\geq 2k+2$, $k\geq 2$ and $p\geq 0$}),   
\end{aligned}
\end{equation*}
which contradicts (\ref{equ::edgetree-2}). It follows that $|V_2|=1$. Let $V_2=\{w\}$. Then $e_{G}(V_1,V_2)=e_{G}(V_1, w)\leq k-1$. Therefore, $G$ is a spanning subgraph of $K_{n,n}\backslash E(K_{1,n-k+1})$. Thus, 
$$\rho(G)\leq \rho(K_{n,n}\backslash E(K_{1,n-k+1})),$$
with equality if and only if $G\cong K_{n,n}\backslash E(K_{1,n-k+1})$, a contradiction.

{\flushleft\bf Case 2.} $t\geq 3$.

In this case, we assert that $|V_1|\geq 2$. Otherwise, $|V_i|=1$ for $1\leq i\leq t$, and hence $t=2n$. Combining this with $n\geq 2k+2$ and (\ref{equ::edgetree-1}), we get 
$$e(G)=e_G(V_1,V_2,\cdots,V_t)\leq k(2n-1)-1\leq \frac{(n-2)(2n-1)}{2}-1<n(n-1),$$
which contradicts (\ref{equ::edgetree-2}). Therefore, $|V_1|\geq 2$, and hence $t\leq 2n-1$.
Suppose that there exist $t_1$ trivial parts in $X$ and $t'_1$ trivial parts in $Y$. Without loss of generality, we may assume that $t_1\geq t_1'$. Thus, $n-t_1\leq n-t'_1$. Let $s$ be the largest index such that $|V_s|\geq 2$, and let $t_2=t'_1+s-1$. Then $n-\sum_{1\leq i\leq s}a_i=t_1$ and $n-\sum_{1\leq i\leq s}b_i=t'_1$.
Combining this with Lemma \ref{lem::1}, (\ref{equ::edgetree-1}) and $n-t_1\leq n-t'_1$, we get
\begin{equation}\label{equ::edgetree-4}
\begin{aligned}
e(G)&= \sum_{1\leq i\leq s} e(V_i)+e_{G}(V_1,V_2,\cdots, V_t)\\
&\leq \sum_{1\leq i\leq s}a_ib_i+k(t-1)-1\\
     &\leq (n-t_1)(n-t'_1-s+1)+k(t-1)-1\\
     &=(n-t_1)(n-t_2)+k(t-1)-1~~(\mbox{since $t_2=t'_1+s-1$})\\
     &=n^2-(t_1+t_2)n+t_1t_2+k(t-1)-1\\
     &\leq n^2-(t-1)n+\frac{(t-1)^2}{4}+k(t-1)-1~~(\mbox{since $t=t_1+t_2+1$})\\
     &=n(n-1)-\Big((t-2)n-\frac{(t-1)^2}{4}-k(t-1)+1\Big).
\end{aligned}
\end{equation}
Let $f(n)=(t-2)n-\frac{(t-1)^2}{4}-k(t-1)+1$. If $t\geq k+2$, then
\begin{equation*}
\begin{aligned}
f(n)&=(t-2)n-\frac{(t-1)^2}{4}-k(t-1)+1\\
    &= \Big(\frac{(t-2)n}{2}-\frac{(t-1)^2}{4}\Big)+\Big(\frac{(t-2)n}{2}-k(t-1)\Big)+1\\
     &\geq \Big(\frac{(t-2)(t+1)}{4}-\frac{(t-1)^2}{4}\Big)+((t-2)(k+1)-k(t-1))+1\\
     &~~(\mbox{since $t\leq 2n-1$ and $n\geq 2k+2$})\\ 
     &=\frac{t-3}{4}+t-k-1\\
     &>0 ~~(\mbox{since $t\geq 3$ and $t\geq k+2$}).
\end{aligned}
\end{equation*}
If $k\geq t-1$, then 
\begin{equation*}
\begin{aligned}
f(n)&=(t-2)n-\frac{(t-1)^2}{4}-k(t-1)+1\\
    &\geq (t-2)(2k+2)-\frac{(t-1)^2}{4}-k(t-1)+1~~(\mbox{since $n\geq 2k+2$ and $t\geq 3$})\\
    &=k(t-3)+2(t-2)-\frac{(t-1)^2}{4}+1\\
    &\geq \frac{3t^2-6t-1}{4}~~(\mbox{since $k\geq t-1$})\\
    &>0 ~~(\mbox{since $t\geq 3$}).
\end{aligned}
\end{equation*}
Therefore, we can deduce that $f(n)>0$ for $n\geq 2k+2$. Combining this with (\ref{equ::edgetree-4}), we have $e(G)<n(n-1)$, which also leads to a contradiction.

This completes the proof.
\end{proof}

\section{Concluding remarks}\label{sec::remarks}

The following fundamental theorem of Hall provides a sufficient and necessary condition for the existence of a perfect matching in a bipartite graph.

\begin{theorem}[Hall \cite{Hall35}]\label{Hall}
An $(X, Y)$-bipartite graph $G$  has a perfect matching if and only if $|X|=|Y|$ and
$$|N_{G}(S)|\geq|S|$$
for all $S\subseteq X$.
\end{theorem}
Motivated by Theorem \ref{Hall}, Fan, Goryainov, Huang and Lin\cite{FGHL22} gave a spectral radius condition for the existence of a perfect matching in a balanced bipartite graph  with fixed minimum degree. 

In fact, Theorem~\ref{Hall} can be extended to matching number.
\begin{theorem}[Ore~\cite{Ore55}]
\label{ore-matching}
The matching number of an $(X, Y)$-bipartite graph $G$ is
$$\alpha'(G) = |X| + \min_{S\subseteq X} \big\{|N(S)|-|S|\big\}.$$
\end{theorem}
Naturally, it would be interesting to study the matching number of bipartite graphs via spectral radius.

%Naturally, we have the following problem.
%\begin{problem}
%Study matching number and spectral radius of bipartite graphs.
%\end{problem}

\section*{Acknowledgements}
Dandan Fan was supported by the National Natural Science Foundation of China (Grant No. 12301454) and sponsored by Natural Science Foundation of Xinjiang Uygur Autonomous Region(Grant Nos. 2022D01B103), Xiaofeng Gu was supported by a grant from the Simons Foundation (522728), and Huiqiu Lin was supported by the National Natural Science Foundation of China (Grant No. 12271162), Natural Science Foundation of Shanghai (No. 22ZR1416300) and The Program for Professor of Special Appointment (Eastern Scholar) at Shanghai Institutions of Higher Learning (No. TP2022031).

\end{document}